\documentclass{article}

\PassOptionsToPackage{square,sort,comma}{natbib}

\usepackage[preprint]{nips_2018}

\usepackage[utf8]{inputenc} 
\usepackage[T1]{fontenc}    
\usepackage[colorlinks = true]{hyperref}       
\usepackage{url}            
\usepackage{booktabs}       
\usepackage{amsfonts}       
\usepackage{nicefrac}       
\usepackage{microtype}      

\usepackage{epsfig,amsmath,amssymb,amsfonts,amstext,amsthm,mathrsfs,dsfont}
\usepackage{latexsym,graphics,epsf,epsfig,psfrag}
\usepackage{color}
\usepackage{subfigure} 
\usepackage{caption}
\usepackage{graphicx}
\usepackage{epstopdf}
\usepackage{algorithm}
\usepackage{algorithmic} 
\usepackage{enumitem} 
\usepackage{cleveref} 
\hypersetup{linkcolor =red, citecolor = blue}

\usepackage{float}

\newcommand{\x}{\mathbf{x}}
\newcommand{\s}{\mathbf{s}}
\newcommand{\ub}{\mathbf{u}}

\newcommand{\norml}[1]{\| #1 \|}
\newcommand{\normlarge}[1]{\left\Vert #1\right\Vert}

\newcommand{\numleqslant}[1]{\overset{\text{(#1)}}{\leqslant}}
\newcommand{\numequ}[1]{\overset{\text{(#1)}}{=}}

\newtheorem{Theorem}{Theorem} 
 
\newtheorem{Lemma}[Theorem]{Lemma} 
\newtheorem{Assumption}{Assumption}

\newcommand\numberthis{\addtocounter{equation}{1}\tag{\theequation}}

\DeclareMathOperator*{\argmin}{argmin}
 
\allowdisplaybreaks

\title{A Note on Inexact Condition for Cubic Regularized Newton's Method}

\author{
	 Zhe Wang\\
	Ohio State University\\	
	wang.10982@osu.edu \\ 
	\And
	Yi Zhou \\
	 Ohio State University\\ 
	 zhou.1172@osu.edu \\
	 \AND
	 Yingbin Liang \\
	 Ohio State University\\
	 liang.889@osu.edu \\
	\And
	Guanghui Lan \\
	Georgia Institute of Technology\\
	george.lan@isye.gatech.edu \\
}
\begin{document}

\maketitle

\begin{abstract}
This note considers the inexact cubic-regularized Newton's method (CR), which has been shown in \cite{Cartis2011a} to achieve the same order-level convergence rate to a secondary stationary point as the exact CR \citep{Nesterov2006}. However, the inexactness condition in \cite{Cartis2011a} is not implementable due to its dependence on future iterates variable. This note fixes such an issue by proving the same convergence rate for nonconvex optimization under an inexact adaptive condition that depends on only the current iterate. Our proof controls the sufficient decrease of the function value over the total iterations rather than each iteration as used in the previous studies, which can be of independent interest in other contexts. 
\end{abstract}
\section{Introduction}
The cubic-regularized (CR) Newton's method \citep{Nesterov2006} is a popular approach that solves the following general nonconvex optimization problem
\begin{align} \label{problem}
  \min_{\x \in \mathbb{R}^d} f(\x),
\end{align}
where $f$ is a differentiable and nonconvex function. Starting from an arbitrary initial point $\x_0$, the update rule of CR can be written as
\begin{align}
	\text{(CR):} \quad \s_{k+1} &= \argmin_{\s \in \mathbb{R}^d} \nabla f(\x_k)^\top\s     + \frac{1}{2} \s^\top \nabla^2 f(\x_k)\s + \frac{M}{6} \norml{\s}^3, \nonumber\\
	\x_{k+1}  &= \x_{k} + \s_{k+1}.
\end{align} 
\cite{Nesterov2006} showed that CR converges to a second-order stationary point $\x$ of the objective function, i.e.,
\begin{align}
	\nabla f(\x) = 0 \quad \text{and} \quad  \nabla^2 f(\x) \succcurlyeq 0.
\end{align}
Such a desirable property allows CR to escape strict saddle points. However, the algorithm  needs to compute a full Hessian at each iteration, and is hence computationally intensive. \cite{Cartis2011a,Cartis2011b} proposed to use an inexact approximation $\mathbf{H}_k$ to replace the full Hessian $\nabla^2 f(\x_{k})$ in the CR update, leading to the following inexact CR algorithm  
\begin{align} 
\hspace{-20mm} \text{(Inexact CR):}   \quad \s_{k+1}  &= \argmin_{\s \in \mathbb{R}^d} \nabla f(\x_k)^\top\s     + \frac{1}{2} \s^\top \mathbf{H}_k \s + \frac{M}{6} \norml{\s}^3, \\
   \x_{k+1}  &= \x_{k} + \s_{k+1}. \label{inexact_cubic_Ca}
\end{align} 
\cite{Cartis2011a,Cartis2011b} showed that if $\mathbf{H}_k$ satisfies the following inexactness condition,
\begin{align}
	\norml{\mathbf{H}_k - \nabla^2 f(\x_{k})} \leqslant C\norml{\s_{k+1}}, \label{future_assumption}
\end{align}
then inexact CR achieves the same order-level convergence rate to a second-order stationary point as exact CR.

This condition has been used in many situations \citep{kohler2017,Cartis2012,Cartis2012b}
Observe that the above inexact condition involves $\norml{\s_{k+1}}$ (and hence $\x_{k+1}$), which is not available at iteration $k$. Thus, the inexact condition in \cref{future_assumption} is not practically implementable. More recent research studies \cite{kohler2017,wang2018} used $s_k$ to replace $s_{k+1}$ to implement inexact CR numerically, i.e., the condition in \cref{future_assumption} becomes
\begin{align}
 	\norml{\mathbf{H}_k - \nabla^2 f(\x_{k})} \leqslant C\norml{\s_{k}}. \label{hisotry_assumption}
\end{align}
These studies demonstrated that inexact CR performs well in experiments under the condition in \cref{hisotry_assumption}, but did not provide theoretical convergence guarantee of inexact CR under such a condition. The main contribution of this note is to establish convergence guarantee for the inexact CR under the inexact condition for Hessian (\cref{hisotry_assumption}) and a similar inexact condition for gradient (see \cref{Inexact_gradient_condition} below), which achieves the same order of convergence rate as the exact CR. In contrast to existing proof techniques, our proof relies on an idea of the overall control of the sufficient decrease of the function value rather than requiring a sufficient decrease at each iteration. More specifically, the inexact error $\norml{\mathbf{H}_k - \nabla^2 f(\x_{k})} \leqslant C\norml{\s_{k}}$ at current iteration is incorporated into the bound on the previous iteration, which yields a successful analysis over all iterations under a more relaxed (and practical) condition \cref{hisotry_assumption}.

We also note that another inexact condition proposed by \cite{Xu2017}, which takes the form 
\begin{align} 
	\norml{\mathbf{H}_k - \nabla^2 f(\x_{k})} \leqslant C \epsilon, \label{conservative_assumption}
\end{align}
where $\epsilon$ is a pre-defined small constant and is required to smaller than  the lower bound of  the set $  \{\norml{\s_i}\}_{i=1}^{k}$  before the algorithm terminates. Similar condition has been used in \cite{Saeed2017,Jiang2017,JinChi2017cubic,Yao2018}. Compared to \cref{conservative_assumption}, \cref{hisotry_assumption} is adapted to $\norml{\s_{k}}$ so that the increment $\norml{\s_{k}}$ can be large than $\epsilon$ in the most phase of the algorithm to enable more progress towards the convergent point.

%

\textbf{Notation:} For a vector $\x$, $\norml{\x}$ denotes the $\ell_2$ norm of the vector $\x$. For a matrix $\mathbf{H}$, $\norml{\mathbf{H}}$ denotes the  spectral norm of the matrix $\mathbf{H}$.  We let $\mathbf{I}$ denote the identity matrix. For a function $f(\cdot)$, $\nabla f(\cdot)$  and $\nabla^2 f(\cdot)$  denote its gradient and Hessian, respectively. $\mathbb{R}$, $\mathbb{R}^+$ and $\mathbb{R}^d$ denote the set of all real numbers,  non-negative real numbers and $d$-dimension real vectors, respectively. $\mathbb{S}$ denotes the set of all symmetric matrices. 
\section{Main Result} 

Our analysis takes the following standard assumption as in the previous studies of CR.
\begin{Assumption} \label{basic_assumption}
	The objective function in \cref{problem} satisfies: 
	\begin{enumerate}[leftmargin=*,topsep=0pt,noitemsep]
		\item   $f(\cdot)$ is twice-continuously differentiable and bounded below, i.e., $f^{\star} \triangleq \inf\limits_{\x \in \mathds{R}^d} f(\x) > -\infty$; 
		\item  The   Hessian $\nabla^2 f(\cdot)$ is $L$-Lipschitz continuous.
	\end{enumerate} 
\end{Assumption}

In our analysis, we allow both the gradient and the Hessian to be replaced by their inexact approximation, and hence the CR iterate becomes
\begin{align} 
\hspace{-20mm}\text{(Inexact gradient and Hessian CR):} \quad \s_{k+1} &= \argmin_{\s \in \mathbb{R}^d} \mathbf{g}_k^\top\s     + \frac{1}{2} \s^\top \mathbf{H}_k \s + \frac{M}{6} \norml{\s}^3, \\
   \x_{k+1}  &= \x_{k} + \s_{k+1}. \label{inexact_cubic_i}
\end{align} 

We assume that $\mathbf{g}_k$ and $\mathbf{H}_k$ satisfy the following inexact conditions, which depend on only the current iteration information, and are hence implementable.
\begin{Assumption}\label{assumption}
There exist two constants $\alpha, \beta \in \mathbb{R}^+$, such that the inexact gradient  $\mathbf{g}_k$ and  inexact Hessian $\mathbf{H}_k$  satisfy, for all $k \geqslant 0$, 
	\begin{align}
	\norml{\mathbf{H}_k - \nabla^2 f(\x_{k})} \leqslant \alpha \norml{\s_{k}}, \\
	\norml{\mathbf{g}_k - \nabla  f(\x_{k})} \leqslant \beta \norml{\s_{k}}^2. \label{Inexact_gradient_condition}
	\end{align} 
\end{Assumption}

%
%

We next state our main theorem, which guarantees that after $k$ iterations, the inexact gradient and Hessian CR must pass an approximate second order saddle point  with error within   $O(1/k^{2/3})$ and $O(1/k^{1/3})$  for the gradient and Hessian, respectively, under the inexact condition in \Cref{assumption}.
\begin{Theorem} \label{Theorem}
	Let Assumptions \ref{basic_assumption} and \ref{assumption} hold. Then, after $k$ iterations, the sequence $\{\x_{i}\}_{i\geqslant 1}$ generated by inexact CR contains a  point $\tilde{\x}$ such that
	\begin{align*}
		 \norml{\nabla f(\tilde{\x})}  \leqslant \frac{C_1}{(k-1)^{2/3}} \quad \text{ and } \quad  \nabla^2 f(\tilde{\x}) \succcurlyeq  -\frac{C_2}{(k-1)^{1/3}} \mathbf{I}.
	\end{align*}
	where $C_1$ and $C_2$ are universal constants, and are specified in the proof.  
\end{Theorem}

The proof of \Cref{Theorem} is based on the following two useful lemmas. 
\begin{Lemma}[\cite{Nesterov2006}, Lemma 1] \label{Hessian_square_bound}
	Let the Hessian $\nabla^2 f(\cdot)$ of the function   $f(\cdot)$ be $L $-Lipschitz continuous with $L  > 0$. Then, for any $\x, \mathbf{y} \in \mathbb{R}^d$, we have
	\begin{align}
	\norml{\nabla f( \mathbf{y} ) - \nabla f(\x) - \nabla^2 f(\x )(\mathbf{y} - \x ) } &\leqslant \frac{L }{2} \norml{\mathbf{y}  - \x}^2.  \label{gradient_bound}\\
	|f(\mathbf{y}) - f(\x) - \nabla f(\x)^T(\mathbf{y}  - \x) - \frac{1}{2} (\mathbf{y}  - \x)^T &\nabla^2 f(\x )(\mathbf{y} - \x )|  \leqslant \frac{L }{6}  \norml{\mathbf{y}  - \x}^3. \label{function_vaule}
	\end{align}
\end{Lemma}

We then establish \Cref{subcubic}, which provides the properties of the minimizer of \cref{inexact_cubic_i} for a more general setting. 
\begin{Lemma} \label{subcubic}
	Let $M \in \mathbb{R}, \mathbf{g} \in \mathbb{R}^d, \mathbf{H} \in  \mathbb{S}^{d \times d}$, and 
	\begin{align}
	\s  = \argmin_{\ub \in \mathbb{R}^d} \mathbf{g}^\top \ub + \frac{1}{2} \ub^\top \mathbf{H} \ub  + \frac{M}{6} \norml{\ub}^3. \label{opt}
	\end{align}
	Then, the following statements  hold:
	\begin{align}
	\mathbf{g}  +      \mathbf{H} \s   + \frac{M}{2} \norml{\s } \s  &= \mathbf{0}, \label{opt_1} \\
	\mathbf{H} + \frac{M}{2} \norml{\s } \mathbf{I}    &\succcurlyeq \mathbf{0}, \label{opt_2} \\
	\mathbf{g}^\top \s  + \frac{1}{2}    \s^\top \mathbf{H} \s   + \frac{M}{6} \norml{\s }^3 &\leqslant - \frac{M}{12}\norml{\s }^3. \label{opt_3}
	\end{align}
\end{Lemma}
To further explain, \cref{opt_1} corresponds to the first-order necessary optimality condition, \cref{opt_2} corresponds to the second-order necessary optimality condition but with a tighter form due to the specific form of this optimization problem, and \cref{opt_3} guarantees a sufficient decrease at this minimizer. 
\begin{proof}[Proof of \Cref{subcubic}]
	First, \cref{opt_1} follows from the first-order necessary optimality condition of \cref{opt}, and \cref{opt_2} follows from the Proposition $1$ in \cite{Nesterov2006}. We next prove \cref{opt_3}. Following similar steps as those in \cite{Nesterov2006}, we obtain that
	\begin{align*}
	\mathbf{g}^\top \s  + \frac{1}{2}   \s^\top \mathbf{H} \s   + \frac{M}{6} \norml{\s }^3 &\numequ{i}    \left(- \mathbf{H} \s   - \frac{M}{2} \norml{\s } \s  \right)^\top \s  + \frac{1}{2}   \s^\top \mathbf{H}  \s      + \frac{M}{6} \norml{\s }^3 \\
	&= - \frac{1}{2}   \s^\top \left(\mathbf{H} +   \frac{M}{2} \norml{\s }\mathbf{I}  \right)\s    - \frac{M}{12}\norml{\s }^3 \numleqslant{ii}  - \frac{M}{12}\norml{\s }^3,
	\end{align*}
	where (i) follows from \cref{opt_1}, and (ii) follows from \cref{opt_2}, which implies that $- \frac{1}{2}   \s^\top \left(\mathbf{H} +   \frac{M}{2} \norml{\s }\mathbf{I}  \right)\s \leqslant 0$.
\end{proof}

Now, we are ready to prove our main theorem.
\begin{proof}[Proof of \Cref{Theorem}]
	Consider any iteration $k$, we obtain that
	\begin{align*}
		f(\x_{k+1}) -f(\x_{k})  &\numleqslant{i}     \nabla f(\x_k)^\top\s_{k+1}     + \frac{1}{2} \s_{k+1}^\top \nabla f(\x_k)\s_{k+1} + \frac{L}{6} \norml{\s_{k+1}}^3 \\
		& \leqslant \mathbf{g}_k^\top\s_{k+1}     + \frac{1}{2} \s_{k+1}^\top \mathbf{H}_k  \s_{k+1}+ \frac{M}{6} \norml{\s_{k+1}}^3 \\
		 &\qquad +  (\nabla f(\x_k) - \mathbf{g}_k)^\top\s_{k+1}  +\frac{1}{2} \s_{k+1}^\top (\nabla^2 f(\x_k) - \mathbf{H}_k)\s_{k+1}  + \frac{L-M}{6} \norml{\s_{k+1}}^3 \\
		&\numleqslant{ii} - \frac{3M-2L}{12}\norml{\s_{k+1} }^3 + (\nabla f(\x_k) - \mathbf{g}_k)^\top\s_{k+1}  + \frac{1}{2} \s_{k+1}^\top (\nabla f(\x_k) - \mathbf{H}_k)\s_{k+1} \\
		&\numleqslant{iii} - \frac{3M-2L}{12}\norml{\s_{k+1} }^3 + \beta \norml{\s_{k}}^2 \norml{\s_{k+1}} + \alpha  \norml{\s_{k}}  \norml{\s_{k+1}}^2 \\
		&\numleqslant{vi} - \frac{3M-2L}{12}\norml{\s_{k+1} }^3 + \beta(\norml{\s_{k}}^3 +\norml{\s_{k+1}}^3 ) + \alpha(\norml{\s_{k}}^3 +\norml{\s_{k+1}}^3 ) \\
		&= - \left( \frac{3M-2L}{12} - \alpha - \beta \right) \norml{\s_{k+1} }^3 + (\alpha + \beta) \norml{\s_{k}}^3. \numberthis \label{iteration}
	\end{align*}  
	where (i) follows from \Cref{Hessian_square_bound} with $\mathbf{y} = \x_{k+1}, \x = \x_{k}$ and $\s_{k+1} =\x_{k+1} - \x_{k}$, (ii) follows from \cref{opt_3} in \Cref{subcubic} with $\mathbf{g} = \mathbf{g}_k,  \mathbf{H} =\mathbf{H}_k $ and $\s = \s_{k+1}$, (iii) follows from \Cref{assumption}, and (vi) follows from the inequality that for $a,b \in \mathbb{R}^+$, $a^2b \leqslant a^3 + b^3$, which can be verified by checking the  cases with $a < b$  and $a \geqslant b$, respectively.
	
	Summing \cref{iteration} from $0$ to $k-1$, we obtain that
	\begin{align*}
			f(\x_{k}) &\leqslant f(\x_{0}) - \sum_{i=0}^{k-1}\left( \frac{3M-2L}{12} - \alpha - \beta \right)\norml{\s_{i+1} }^3  +  \sum_{i=0}^{k-1} (\alpha + \beta) \norml{\s_{i}}^3, \\
			&\leqslant f(\x_{0}) - \sum_{i=1}^{k }\left( \frac{3M-2L}{12} - \alpha - \beta \right)\norml{\s_{i} }^3  +  \sum_{i=0}^{k } (\alpha + \beta) \norml{\s_{i}}^3 \\
			&\leqslant f(\x_{0}) - \sum_{i=1}^{k }\left( \frac{3M-2L}{12} - \alpha - \beta \right)\norml{\s_{i} }^3  +  \sum_{i=1}^{k} (\alpha + \beta) \norml{\s_{i}}^3  +    (\alpha + \beta) \norml{\s_{0}}^3 \\
			&= f(\x_{0}) - \sum_{i=1}^{k }\left(\frac{3M-2L}{12} - 2\alpha - 2\beta \right)\norml{\s_{i} }^3   +    (\alpha + \beta) \norml{\s_{0}}^3.
	\end{align*}
	Therefore, we have \begin{align}
		\sum_{i=1}^{k } \gamma \norml{\s_{i} }^3 \leqslant f(\x_{0}) -	f^\star + (\alpha + \beta) \norml{\s_{0}}^3,  \label{total_sum}
	\end{align}
	where $\gamma \triangleq  \frac{3M-2L}{12} - 2\alpha - 2\beta$. We note that $M > \frac{2}{3}L + 8\alpha + 8\beta$, and thus we have $\gamma = \frac{3M-2L}{12} - 2\alpha - 2\beta > 0$. Let $ m \triangleq \argmin_{i \in \{1,  \cdots, k -1   \}} \norml{\s_{i} }^3 + \norml{\s_{i+1} }^3$. We obtain that
	\begin{align*}
		 \norml{\s_{m} }^3 + \norml{\s_{m+1} }^3 &= \min_{i \in \{1,  \cdots, k -1    \}} \norml{\s_{i} }^3 + \norml{\s_{i+1} }^3 \\
		  &\leqslant \frac{1}{k-1} \sum_{i=1}^{k-1}\left( \norml{\s_{i} }^3  + \norml{\s_{i+1} }^3  \right) \\ &\numleqslant{i} \frac{2}{\gamma (k-1)}  \left({f(\x_{0}) -	f(\x^\star) + (\alpha + \beta) \norml{\s_{0}}^3} \right).
	\end{align*}
	where (i) follows \cref{total_sum}.
	
	Therefore, we have
	\begin{align}
		 \max \left\{\norml{\s_{m} },\norml{\s_{m+1} } \right\} \leqslant \frac{1}{(k-1)^{1/3}}\left( \frac{2}{\gamma} \left({f(\x_{0}) -	f(\x^\star) + (\alpha + \beta) \norml{\s_{0}}^3} \right)\right)^{1/3}.    \label{norm_convergence}
	\end{align}
	
	Next, we prove the convergence rate of $\nabla f(\cdot)$ and $\nabla^2 f(\cdot)$. We first derive 
	\begin{align*}
	  \norml{\nabla f(\x_{m+1})} &\numequ{i} \normlarge{ \nabla f(\x_{m+1}) - \left(\mathbf{g}_m + \mathbf{H}_m \s_{m+1} + \frac{M}{2}\norml{\s_{m+1}} \s_{m+1} \right) } \\
	  &\leqslant \normlarge{ \nabla f(\x_{m+1}) - \left(\mathbf{g}_m + \mathbf{H}_k \s_{m+1} \right) } + \frac{M}{2} \norml{\s_{m+1}}^2 \\
	  &\leqslant \normlarge{ \nabla f(\x_{m+1}) - \nabla f(\x_m) -  \nabla^2 f(\x_m)\s_{m+1} } \\
	     &\qquad + \norml{\nabla f(\x_m) - \mathbf{g}_m } + \norml{(\nabla^2 f(\x_m) - \mathbf{H}_m )\s_{m+1}} + \frac{M}{2} \norml{\s_{m+1}}^2 \\
	  &\numleqslant{ii} \frac{L}{2} \norml{\s_{m+1}}^2 + \beta \norml{\s_{m}}^2 + \alpha \norml{\s_{m}}\norml{\s_{m+1}}  + \frac{M}{2} \norml{\s_{m+1}}^2 \\
	  &\numleqslant{iii}  \frac{1}{(k-1)^{2/3}} \frac{L + M + 2\beta + 2\alpha}{2}  \left( \frac{2}{\gamma} \left({f(\x_{0}) -	f^\star + (\alpha + \beta) \norml{\s_{0}}^3} \right)\right)^{2/3},
	\end{align*}
	where (i) follows from \cref{opt_1} with $\mathbf{g} = \mathbf{g}_m, \mathbf{H} = \mathbf{H}_m$ and $\s = \s_{m+1}$, (ii) follows from \cref{gradient_bound} in \Cref{Hessian_square_bound} and \Cref{assumption}, and (iii) follows from \cref{norm_convergence}.
	
    We next prove the   the convergence rate of $\nabla^2 f(\cdot)$.
    \begin{align*}
     \nabla^2 f(\x_{m+1}) &\overset{(i)}{\succcurlyeq} \mathbf{H}_m - \norml{\mathbf{H}_m -\nabla^2 f(\x_{m+1})} \mathbf{I}  \\
     &\overset{(ii)}{\succcurlyeq} - \frac{M}{2} \norml{\s_{m+1}} \mathbf{I}  - \norml{\mathbf{H}_m -\nabla^2 f(\x_{m+1})} \mathbf{I} \\
     &\succcurlyeq  - \frac{M}{2} \norml{\s_{m+1}} \mathbf{I}  - \norml{\mathbf{H}_m -\nabla^2 f(\x_{m })}\mathbf{I}  - \norml{\nabla^2 f(\x_{m}) -\nabla^2 f(\x_{m+1})}\mathbf{I}  \\
     &\overset{(iii)}{\succcurlyeq}  - \frac{M}{2} \norml{\s_{m+1}} \mathbf{I}  -  \alpha \norml{\s_{m}}\mathbf{I}  -  L \norml{\s_{m+1}}\mathbf{I} \\
     &\overset{(iv)}{\succcurlyeq} -\frac{1}{(k-1)^{1/3}} \frac{M + 2L + 2\alpha}{2}   \left( \frac{2}{\gamma} \left({f(\x_{0}) -	f(\x^\star) + (\alpha + \beta) \norml{\s_{0}}^3} \right)\right)^{1/3} \mathbf{I} ,
    \end{align*}
	where (i) follows from Weyl's inequality, (ii) follows from \cref{opt_2} with $\mathbf{H} = \mathbf{H}_{m}$ and $\s = \s_{m+1}$, (iii) follows from \Cref{assumption} and the fact that $\nabla^2 f(\cdot)$ is $L-$Lipschitz, and (vi) follows from \cref{norm_convergence}. 
\end{proof}

\section{Conclusion}
In this note, we study the cubic-regularized Newton's method under a more practical inexact condition, which depends only on the current iteration information, rather than the future iteration in previous studies. Under such an inexact condition for both the gradient and the Hessian, we establish the convergence of the inexact CR method to a second-order stationary point, and show that the convergence rate is as fast as that of CR in nonconvex optimization.

\bibliographystyle{apalike} 
\bibliography{BibTex}

\begin{thebibliography}{}

\bibitem[Cartis et~al., 2012a]{Cartis2012}
Cartis, C., Gould, N., and Toint, P.~L. (2012a).
\newblock Complexity bounds for second-order optimality in unconstrained
  optimization.
\newblock {\em Journal of Complexity}, 28(1):93 -- 108.

\bibitem[Cartis et~al., 2011a]{Cartis2011a}
Cartis, C., Gould, N. I.~M., and Toint, P.~L. (2011a).
\newblock {Adaptive cubic regularization methods for unconstrained
  optimization. Part I : Motivation, convergence and numerical results}.
\newblock {\em Mathematical Programming}, 127(2):245--295.

\bibitem[Cartis et~al., 2011b]{Cartis2011b}
Cartis, C., Gould, N. I.~M., and Toint, P.~L. (2011b).
\newblock {Adaptive cubic regularization methods for unconstrained
  optimization. Part II worst-case function- and derivative-evaluation
  complexity}.
\newblock {\em Mathematical Programming}, 130(2):295--319.

\bibitem[Cartis et~al., 2012b]{Cartis2012b}
Cartis, C., Gould, N. I.~M., and Toint, P.~L. (2012b).
\newblock An adaptive cubic regularization algorithm for nonconvex optimization
  with convex constraints and its function-evaluation complexity.
\newblock {\em IMA Journal of Numerical Analysis}, 32(4).

\bibitem[Ghadimi et~al., 2017]{Saeed2017}
Ghadimi, S., Liu, H., and Zhang, T. (2017).
\newblock Second-order methods with cubic regularization under inexact
  information.
\newblock {\em arXiv: 1710.05782}.

\bibitem[{Jiang} et~al., 2017]{Jiang2017}
{Jiang}, B., {Lin}, T., and {Zhang}, S. (2017).
\newblock {A unified scheme to accelerate adaptive cubic regularization and
  gradient methods for convex optimization}.
\newblock {\em arXiv:1710.04788}.

\bibitem[Kohler and Lucchi, 2017]{kohler2017}
Kohler, J.~M. and Lucchi, A. (2017).
\newblock Sub-sampled cubic regularization for non-convex optimization.
\newblock In {\em Proc. 34th International Conference on Machine Learning
  (ICML)}, volume~70, pages 1895--1904.

\bibitem[Nesterov and Polyak, 2006]{Nesterov2006}
Nesterov, Y. and Polyak, B.~T. (2006).
\newblock Cubic regularization of newton method and its global performance.
\newblock {\em Mathematical Programming}, 108(1):177--205.

\bibitem[{Tripuraneni} et~al., 2017]{JinChi2017cubic}
{Tripuraneni}, N., {Stern}, M., {Jin}, C., {Regier}, J., and {Jordan}, M.~I.
  (2017).
\newblock Stochastic cubic regularization for fast nonconvex optimization.
\newblock {\em arXiv: 711.02838}.

\bibitem[{Wang} et~al., 2018]{wang2018}
{Wang}, Z., {Zhou}, Y., {Liang}, Y., and {Lan}, G. (2018).
\newblock Sample complexity of stochastic variance-reduced cubic regularization
  for nonconvex optimization.
\newblock {\em arXiv:1802.07372}.

\bibitem[Xu et~al., 2017]{Xu2017}
Xu, P., Roosta-Khorasani, F., and Mahoney, M.~W. (2017).
\newblock Newton-type methods for non-convex optimization under inexact hessian
  information.
\newblock {\em arXiv: 1708.07164}.

\bibitem[{Yao} et~al., 2018]{Yao2018}
{Yao}, Z., {Xu}, P., {Roosta-Khorasani}, F., and {Mahoney}, M.~W. (2018).
\newblock {Inexact Non-Convex Newton-Type Methods}.
\newblock {\em ArXiv:1802.06925}.

\end{thebibliography}

\end{document}